\documentclass[12pt,a4paper,oneside,reqno,notitlepage]{amsart}

\usepackage{amscd,amsmath}

\usepackage{amssymb,latexsym}
\usepackage{eucal}
\usepackage{graphics}

\usepackage[arrow,matrix,curve]{xy}\SilentMatrices
\def\xyma{\xymatrix@M.7em}

\numberwithin{equation}{section}

\newtheorem{cor}{Corollary}[section]

\newtheorem{prop}{Proposition}[section]
\newtheorem{theorem}{Theorem}[section]

\def\bee{\begin{equation}}
\def\ee{\end{equation}}

\begin{document}

\title{The subgroup determined by a certain ideal in a free group ring}
\begin{quote}
\begin{abstract}
For normal subgroups   $R$ and $S$ of a free group $F$,
an\linebreak identification of the subgroup $F\cap (1+\mathfrak
r\mathfrak f\mathfrak s)$ is derived, and it is shown that the the
quotient  $\frac{F\cap (1+\mathfrak r\mathfrak f\mathfrak
s)}{[R'\cap S',\, R\cap S][R'\cap S,\,R'\cap S][R\cap S',\, R\cap
S']}$ is, in general, non-trivial.
\end{abstract}\end{quote}
\author{Roman Mikhailov and Inder Bir S. Passi}
\maketitle
\par\vspace{.25cm}\noindent
Key words and phrases: Group ring, augmentation ideal, free group, derived functors of non-additive functors, \par\vspace{.25cm}\noindent MSC2010:  18E25, 20C05, 20C07, 20E05.

\section{Introduction}
Every two-sided ideal $\mathfrak a$ in the integral group ring
$\mathbb Z[F] $ of a free group $F$ determines a normal subgroup
$F \cap (1 + \mathfrak a)$ of $F$. Identification of such
subgroups is a fundamental problem in the theory of group rings
(\cite{Gupta:87}, \cite{MP:2009}). Let $R$ and $S$ be normal
subgroups of $F$. In this paper we examine   the subgroup $F\cap
(1+\mathfrak r\mathfrak f\mathfrak s)$, where, for a normal
subgroup $G$ of $F$, $\mathfrak g$ denotes the two-sided  ideal of
$\mathbb Z[F]$ generated by $G-1$. This subgroup has been studied
by \linebreak C. K. Gupta \cite{CKGupta} (see also
\cite{KK:2002}). It is easy to check that $$[R'\cap S',\, R\cap
S][R'\cap S,\,R'\cap S][R\cap S',\, R\cap S']\subseteq F\cap
(1+\mathfrak r\mathfrak f\mathfrak s),$$ where $R'$ (resp. $S'$)
is the derived subgroup of $R$ (resp. $S$).   Whereas the
identification given in \cite{CKGupta}, namely that the preceding
inclusion is an equality,  holds up to torsion,  our investigation
shows that, $\frac{F\cap (1+\mathfrak r\mathfrak f\mathfrak
s)}{[R'\cap S',\, R\cap S][R'\cap S,\,R'\cap S][R\cap S',\, R\cap
S']}\cong L_1{\sf SP}^2\left(\frac{R\cap S}{(R'\cap S)(R\cap
S')}\right)$, and is, in general,  non-identity;   here $L_1{\sf
SP}^2$ is the first derived functor of the second symmetric power
functor.
\par\vspace{.5cm}

 \section{The subgroup $F\cap (1+\mathfrak r\mathfrak f\mathfrak s)$}\par\vspace{.5cm}
 Let $F$ be a free group and $R,\ S$ its normal subgroups with bases, as free groups,  $\{r_i\}_{i\in I}$ and $\{s_j\}_{j\in J}$ respectively.  Then the ideal $\mathfrak r$ is a free right $\mathbb Z[F]$-module with\linebreak  basis $\{r_i-1\,|\,i\in I\}$ and the ideal $\mathfrak s$  is a free left  $\mathbb Z[F]$-module with basis $\{s_j-1\,|\,j\in J\}$ \linebreak (\cite{Gruenberg:1970}, Theorem 1, p.\,32).  Further, recall that  $$
R/R'\cong \frac{\mathfrak  r}{\mathfrak r\mathfrak f}\ \ \ (\cong \frac{\mathfrak  r}{\mathfrak
f\mathfrak r})
$$ this isomorphism being given by $$rR'\mapsto (r-1)+{\mathfrak r\mathfrak f},\ r\in R\quad (\text{resp.} \ rR'\mapsto r-1+\mathfrak f\mathfrak r).$$
From these observations it immediately follows that we can make the following \linebreak identification
\begin{equation}\label{tensor} R/R'\otimes S/S'=\frac{\mathfrak  r}{\mathfrak r\mathfrak f}\otimes
\frac{\mathfrak  s}{\mathfrak f\mathfrak s}=\frac{\mathfrak  r}{\mathfrak r\mathfrak f}\otimes_{\mathbb
Z[F]}\frac{\mathfrak  s}{\mathfrak f\mathfrak s}=\frac{\mathfrak  r\mathfrak s}{\mathfrak r\mathfrak f\mathfrak s}.
\end{equation}
Here $\otimes$ is tensor product over $\mathbb Z$ which we can replace by  $\otimes_{\mathbb Z[F]}$ since the
action of $\mathbb Z[F]$ on components of the tensor product is
trivial.

\par\vspace{.5cm}\noindent
\begin{theorem}\label{the1} If  $R$ and  $S$ are  normal subgroups of a free group $F$, then
there is a natural isomorphism $$ \frac{F\cap (1+{\mathfrak
r\mathfrak f\mathfrak s})}{[R'\cap S', R\cap S][R'\cap S,R'\cap
S][R\cap S', R\cap S']}\cong L_1{\sf SP}^2\left(\frac{R\cap
S}{(R'\cap S)(R\cap S')}\right).
$$
\end{theorem}

\par\vspace{.5cm}\noindent
{\bf Proof.}
Let us set \begin{equation}\label{notation1} Q:=\frac{R\cap
S}{R'\cap S'},\quad U:=\frac{R'\cap S}{R'\cap S'},\quad
V:=\frac{R\cap S'}{R'\cap S'}.\end{equation} The group $Q$ is
free abelian because it injects into $R/R'\oplus S/S'$, and so are
$U,\, V$ both being subgroups of $Q$. Observe that $Q/U$ is also
free abelian, since it is isomorphic to the subgroup $\frac{R\cap
S}{R'\cap S}$ of $R/R'$.
\par\vspace{.25cm}
For an abelian group $A$, we denote by ${\sf SP}^2(A)$ its
symmetric square, defined as the quotient ${\sf SP}^2(A):=A\otimes
A/\langle a\otimes b-b\otimes a,\,|\, a,\,b\in A\rangle$ and by $\Lambda^2(A)$ its
exterior square $\Lambda^2(A):=A\otimes A/\langle a\otimes a\,|\, a\in A\rangle$.
Recall  (see \cite{Kock:2001}) that, for any free resolution
 $$0\to C\to B\to A\to 0$$ of $A$, the so-called {\it Koszul complex}
$$
0\to \Lambda^2(C)\to C\otimes B\to {\sf SP}^2(B)
$$
represents the object $L{\sf SP}^2(A)$ of the derived category of
abelian groups; in particular, its zeroth (resp. first) homology is
equal to the zeroth (resp. first) derived functor of ${\sf SP}^2$
applied to $A$.
\par\vspace{.25cm}

Consider the natural commutative diagram with exact rows and
columns which contains maps between quadratic Koszul complexes:
$$
\xyma{\Lambda^2(U)\ar@{>->}[d] \ar@{>->}[r] & U\otimes
Q\ar@{->}[r] \ar@{>->}[d] & {\sf SP}^2(Q)\ar@{=}[d]\\
\Lambda^2(Q)\ar@{->}[r]\ar@{->>}[d] & Q\otimes Q\ar@{->}[r] \ar@{->>}[d] & {\sf SP}^2(Q)\\
\frac{\Lambda^2(Q)}{\Lambda^2(U)}\ar@{->}[r] & Q/U\otimes Q}
$$
Since the middle horizontal complex is acyclic, the homology of
the lower complex are the same as of the upper complex shifted by
one. That is, there exists a short exact sequence
$$
0\to \frac{\Lambda^2(Q)}{\Lambda^2(U)}\to Q/U\otimes Q\to {\sf
SP}^2(Q/U)\to 0
$$
which can be naturally extended to the following diagram:
\begin{equation}\label{di1}
\xyma{& & & K\ar@{>->}[d]\\ & & Q/U\otimes V\ar@{=}[r] \ar@{>->}[d] & Q/U\otimes V\ar@{->}[d]\\
& \frac{\Lambda^2(Q)}{\Lambda^2(U)}\ar@{>->}[r] \ar@{=}[d] &
Q/U\otimes
Q\ar@{->>}[r] \ar@{->>}[d] & {\sf SP}^2(Q/U)\\
K\ar@{>->}[r] & \frac{\Lambda^2(Q)}{\Lambda^2(U)} \ar@{->}[r] &
Q/U\otimes Q/V}
\end{equation}
Here $K$ is, by definition, the kernel of the lower horizontal
map. By Snake Lemma, $K$ is isomorphic to the kernel of the right
hand vertical map $Q/U\otimes V\to {\sf SP}^2(Q/U) $ in the diagram. Observe that this map is  part of
the Koszul complex
$$
0\to \Lambda^2(VU/U)\to Q/U\otimes V\to {\sf SP}^2(Q/U)
$$
which represents the object $L{\sf SP}^2(Q/UV)$ of the derived
category of abelian groups. Here we have used the fact that
$V=VU/U=V/(V\cap U),$ since $V\cap U$ is the zero subgroup of $Q$.
The homology groups of the above Koszul complex are the derived functor evaluations
$L_i{\sf SP}^2(Q/UV),\ i=1,\,2$ (see \cite{Kock:2001}). Therefore, we get the following
short exact sequence:
$$
0\to \Lambda^2(V)\to K\to L_1{\sf SP}^2(Q/UV)\to 0.
$$
Consequently  the lower sequence of the diagram (\ref{di1}),
yields  the following exact sequence: \begin{equation}\label{qw}
0\to L_1{\sf SP}^2(Q/UV)\to
\frac{\Lambda^2(Q)}{\Lambda^2(U)+\Lambda^2(V)}\to Q/U\otimes Q/V
\end{equation}
We next observe that there are natural isomorphisms
\begin{align*}
& \Lambda^2(Q)\cong \frac{\gamma_2(R\cap S)}{[R'\cap S',\,R\cap S]}\\
& \frac{\Lambda^2(Q)}{\Lambda^2(U)+\Lambda^2(V)}\cong
\frac{\gamma_2(R\cap S)}{[R'\cap S',\,R\cap S][R'\cap S,\,R'\cap
S][R\cap S',\,R\cap S']},
\end{align*}
and natural monomorphisms $Q/U\to R/R'$, $Q/V\to S/S'$. The exact
sequence (\ref{qw}) thus implies that there is an exact sequence
\begin{equation}\label{mn1}
0\to L_1{\sf SP}^2(Q/UV)\to \frac{\gamma_2(R\cap S)}{[R'\cap
S',R\cap S][R'\cap S,R'\cap S][R\cap S',R\cap S']}\to R/R'\otimes
S/S'.
\end{equation}
The statement of the theorem follows from the fact (see \cite{BD})
that
$$
F\cap (1+{\mathfrak  r\mathfrak s})=\gamma_2(R\cap S)
$$
and the identification (\ref{tensor}).\ \ $\Box$
\par\vspace{.5cm}

For an abelian group $A$, a description of the group $L_1{\sf
SP}^2(A)$ is available in many papers on polynomial functors; for
example, see \cite{BP} or (\cite{Jean:2002}, Theorem 2.2.5).
Recall the main properties of $L_1{\sf SP}^2(A)$. For any abelian group  $A$,
$L_1{\sf SP}^2(A)$ is a natural quotient of the group ${\sf
Tor}(A,\,A)$ by diagonal elements. We have
$$L_1{\sf SP}^2(\mathbb Z/m\mathbb Z)=L_1{\sf SP}^2(\mathbb Z)=0,\
$$ for all natural numbers $m$,  and, for all abelian groups $A,\,B$, there is a
(bi)natural isomorphism
$$
{\sf Tor}(A,\,B)={\sf Ker}\{L_1{\sf SP}^2(A\oplus
B)\twoheadrightarrow L_1{\sf SP}^2(A)\oplus L_1{\sf SP}^2(B)\}.
$$
For a free abelian group $A$ and a natural number $m\geq 1$, there is a natural
isomorphism
$$
L_1{\sf SP}^2(A\otimes \mathbb Z/m\mathbb Z)\simeq
\Lambda^2(A\otimes\mathbb Z/m\mathbb Z).
$$
Observe also that, the functor $L_1{\sf SP}^2$ is related to the
homology of the Eilenberg-MacLane spaces $K(-,\,2)$. Namely, for any
abelian group $A$, there is a natural short exact sequence
$$
0\to L_1{\sf SP}^2(A)\to H_5K(A,\,2)\to {\sf Tor}(A,\,\mathbb Z/2\mathbb Z)\to
0.
$$

Invoking this description for $L_1{\sf SP}^2(Q/UV)$, we have the
following identification of the subgroup $F\cap (1+\mathfrak
r\mathfrak f\mathfrak s)$:\par\vspace{.5cm}
\begin{theorem}
$$F\cap (1+{\mathfrak  r\mathfrak f\mathfrak s})=[R'\cap S',\, R\cap S][R'\cap S,\,R'\cap
S][R\cap S', \,R\cap S']W,$$ where $W$ is the subgroup of $F$
generated by elements\footnote{For an elements $g,\,h$ of a group,
we use the standard commutator notation
$[g,\,h]:=g^{-1}h^{-1}gh$.}
$$[x_1,\,y][x,\,y_2]^{-1},$$ such
that \begin{align*} & x,\,y\in R\cap S,\ m\geq 2,\\
& x^m=x_1x_2,\ y^m=y_1y_2,\\ & x_1,\,y_1\in R'\cap S,\\ & x_2,\,y_2\in
R\cap S'.\end{align*}
\end{theorem}
\par\vspace{.25cm}\noindent{\bf Proof.}
Consider the generating elements from $W$, as in the Theorem. Modulo
$\mathfrak r\mathfrak f\mathfrak s$, we have
\begin{align*}
[x_1,\,y][x,\,y_2]^{-1}-1& \equiv [x^m,\,y][x_2,\,y]^{-1}[x,\,y_2]^{-1}-1\\
& \equiv (x^m-1)(y-1)-(y^m-1)(x-1)\\ & \ \ \ \ -(x_2-1)(y-1)+(y-1)(x_2-1)\\
& \ \ \ \ -(x-1)(y_2-1)+(y_2-1)(x-1)\\
& \equiv (x_1-1)(y-1)-(y_1-1)(x-1)+\\ & \ \ \ \
(y-1)(x_2-1)-(x-1)(y_2-1).
\end{align*}
All four products
$(x_1-1)(y-1),\,(y_1-1)(x-1),\,(y-1)(x_2-1),\,(x-1)(y_2-1)$ lie in
$\mathfrak r\mathfrak f\mathfrak s$. The subgroup $W$ is chosen as
a subgroup of representatives of $L_1{\sf SP}^2\left(\frac{R\cap
S}{(R'\cap S)(R\cap S')}\right)$ in $F\cap (1+\mathfrak r\mathfrak
f\mathfrak s)$.
\par\vspace{.25cm}
Consider generators of $L_1{\sf SP}^2\left(\frac{R\cap S}{(R'\cap
S)(R\cap S')}\right)$ viewed as a natural quotient of the group
${\sf Tor}\left(\frac{R\cap S}{(R'\cap S)(R\cap S')},\, \frac{R\cap
S}{(R'\cap S)(R\cap S')}\right).$ The generators are given as
pairs of elements $(x,\,y),\ x,\,y\in R\cap S,$ with the property that,
there exists $m\geq 2$, such that $x^m,\, y^m \in (R'\cap S)(R\cap
S')$. Consider now the diagram (\ref{di1}) and find the image of
the pair $(x,\,y)$ in the quotient $\frac{Q/U\otimes
Q}{\Lambda^2(V)}$ (here we use the notation \ref{notation1}) and
choose its representative in $Q/U\otimes V$. It is given as
\begin{equation}\label{repres}(x.U)\otimes y_2.(R'\cap
S')-(y.U)\otimes x_2.(R'\cap S'),\end{equation} where $x_2,\,y_2$
are defined in the formulation of the Theorem. Going further in
the diagram (\ref{di1}), we find a representative of the element
(\ref{repres}) in $\Lambda^2(Q)/\Lambda^2(V)$, given as
$$
(x\wedge y_2)+(x_2\wedge y)-(x^m\wedge y)+\Lambda^2(V).
$$
Indeed, the natural map $\Lambda^2(Q)\to Q/U\otimes Q$ sends (we
omit the notation $-.(R'\cap S')$ for the elements from $Q$ for the sake of simplification of  notations)
\begin{align*}
(x\wedge y_2)+(x_2\wedge y)-(x^m\wedge y)\mapsto\ \ & x.U\otimes
y_2-y_2.U\otimes x+x_2.U\otimes y-y.U\otimes x_2\\
& -x_1x_2.U\otimes
y-y.U\otimes x_1x_2=\\
& x.U\otimes y_2-y_2.U\otimes x-y.U\otimes x_1=\\
& x.U\otimes y_2-y^m.U\otimes x-y.U\otimes x_1=\\
& x.U\otimes y_2-y.U\otimes x^m-y.U\otimes x_1=\\
& x.U\otimes y_2-y.U\otimes x_2.
\end{align*}
In the free group $F$, this element is represented
as a product of commutators $$[x,\,y_2][x_2,\,y][x^m,\,y]^{-1}.$$ Since
modulo $\mathfrak r\mathfrak f\mathfrak s,$
$$[x_1,\,y][x,\,y_2]^{-1}-1\equiv
[x^m,y][x_2,\,y]^{-1}[x,\,y_2]^{-1}-1\equiv
([x,\,y_2][x_2,\,y][x^m,\,y]^{-1})^{-1}-1$$ we get the asserted
description of the set $W$. $\Box$

\vspace{.5cm}\noindent{\bf Remark.} Since the groups
$F/\gamma_2(R\cap S),\  R/R', \ S/S'$ are always torsion-free, the
sequence (\ref{mn1}) implies that there is the following
identification
$$
L_1{\sf SP}^2\left(\frac{R\cap S}{(R'\cap S)(R\cap
S')}\right)\cong {\sf torsion\ of}\ \frac{F}{[R'\cap S',R\cap
S][R'\cap S,R'\cap S][R\cap S',R\cap S']}.
$$

\par\vspace{.5cm}
\section{Example}
\par\vspace{.5cm} Finally, let us give an example of subgroups $R,\ S$ in a
free group $F$, such that
$$
L_1{\sf SP}^2\left(\frac{R\cap S}{(R'\cap S)(R\cap S')}\right)\neq
0.
$$
Let $F=F(a_1,\,\dots\,,\, a_n,\,b)$, $n\geq 2$,
\begin{align*}
& R=\langle a_1,\,\dots\,,\, a_n,\, [F,\,F]\rangle^F,\\
& S=\langle a_1^2,\,\dots\,,\, a_n^2,\, b,\,[F,\,F]\rangle^F.
\end{align*}
Since $[F,\,F]\subset R,\, [F,\,F]\subset S,$
$$
(R'\cap S)(R\cap S')=R'S'.
$$
For every $i=1,\,\dots\,,\, n$, the element $[a_i,\,b]$ lies in $R\cap S$.
Observe that,
$$
[a_i^2,\,b]=[a_i,\,b][[a_i,\,b],\,a_i][a_i,\,b]
$$
Therefore,
$$
[a_i,\,b]^2\in R'S'.
$$
Since $R'S'=\langle [a_i,\,a_j], \,[a_i,\,b]^2, \gamma_3(F)\rangle$, the
elements $[a_i,\,b],\ i=1,\,\dots\,,\,n$ form an abelian subgroup of
$\frac{R\cap S}{(R'\cap S)(R\cap S')}$ isomorphic to $(\mathbb
Z/2)^{\oplus n}$. For $n\geq 2$, the first derived functor of
${\sf SP}^2$ of such group is non-zero.

\section*{Acknowledgement}
The research of the first author is supported by the Russian
Science Foundation, grant N 14-21-00035. The authors thank S. O.
Ivanov for discussions related to the subject of the paper.

\par\vspace{.5cm}\noindent
Roman Mikhailov\\
St Petersburg Department of Steklov Mathematical Institute\\
and\\
Chebyshev Laboratory\\
St Petersburg State University\\
14th Line, 29b\\
Saint Petersburg\\
199178 Russia\\
email: romanvm@mi.ras.ru
\par\vspace{.5cm}\noindent
Inder Bir S. Passi\\
Centre for Advanced Study in Mathematics\\
Panjab University\\
Sector 14\\
Chandigarh 160014 India\\
and\\
Indian Institute of Science Education and Research\\
Mohali (Punjab)140306 India\\
email: ibspassi@yahoo.co.in

\end{document}